\documentclass{article}
\usepackage{amsfonts}
\usepackage{amsmath}

\input{tcilatex}

\begin{document}

\title{Three competing patterns}
\date{}
\author{Rita Abraham and Jan \ Vrbik \\
Department of Mathematics\\
Brock University, 500 Glenridge Ave.\\
St. Catharines, Ontario, Canada, L2S 3A1}
\maketitle

\begin{abstract}
Assuming repeated independent sampling from a Bernoulli distribution with
two possible outcomes $S$ and $F,$ there are formulas for computing the
probability of one specific pattern of consecutive outcomes (such as $SSFFSS$%
) winning (i.e. being generated first) over another such pattern (e.g. $%
SFSSFS$). In this article we will extend the theory to three competing
patterns.
\end{abstract}

\date{}

\section{Completing a pattern}

\subsection{from scratch}

Consider a sequence of \emph{independent} Bernoulli trials, each having two
potential outcomes, \textsf{S} (a success) with the probability of $p$ and 
\textsf{F} (a failure) with the probability of $q=1-p.$ Then, choose a
specific pattern of $m$ consecutive outcomes (such as \textsf{SSFFS}), and
define the \textsc{probability generating function} (PGF for short) of the
number of trials needed to generate this pattern \emph{for the first time,
from scratch} (let us call the corresponding random variable $X$) as%
\begin{equation*}
F(s)\equiv \mathbb{E}\left( s^{X}\right) =\sum_{n=0}^{\infty }f_{n}s^{n}
\end{equation*}%
where $f_{n}$ is the probability of completing the pattern, for the first
time, at Trial $n.$

If $u_{n}$ is the probability of completing the same pattern (but not
necessarily for the first time) at Trial $n,$ we can relate the two
sequences by%
\begin{equation}
u_{n}=\sum_{i=0}^{n}f_{i}u_{n-i}  \label{UF}
\end{equation}%
which is correct for any $n\geq 1$ (but not for $n=0$) after setting $u_{0}=1
$ and $f_{0}=0$ (to prove that, partition the sample space of the first $n$
trials according the trial - denoted $i$ - at which the first occurrence of
the pattern happened, and use the total-probability formula). Here we are
assuming (and this is quite crucial) that upon completing one occurrence of
the pattern, we are \emph{not allowed} to use any of its symbols to help
build its next occurrence - we have to start `from scratch'. Let $U(s)$ be
the sequence generating function (SGF) of the $u_{n}$ probabilities, i.e. $%
U(s)\equiv \sum_{n=0}^{\infty }u_{n}s^{n}$; note that $U(1)=\infty $ - these
probabilities do not constitute a distribution the way the $f_{n}$
probabilities do.

Multiplying (\ref{UF}) by $s^{n}$ and summing over $n$ from $1$ to infinity
yields $U(s)-1$ (to account for the missing $u_{0}$) on the LHS and $%
U(s)\cdot F(s)$, representing the \textsc{convolution} of the two sequences,
on the RHS (simplifying a convolution of two sequences in this manner is the
main reason for using generating functions in this context). The resulting
equation can be easily solved for $F(s),$ leading to%
\begin{equation}
F(s)=\frac{U(s)-1}{U(s)}  \label{FU}
\end{equation}

To utilize this formula, we must first find $U(s),$ which happens to be
easier than finding, directly, $F(s)$. To achieve the former, we relate the
probability of finding the \emph{symbols} of the pattern (visualize \textsf{%
SSFFS)} at Trials $n-m+1$ to $n$ (note that this does \emph{not} necessarily
mean that the corresponding pattern has been completed at Trial $n$ - it may
have been completed \emph{earlier}, e.g. at Trial $n-4$ using our \textsf{%
SSFFS} example, which would prevent its completion at Trial $n,$ due to the
from-scratch requirement) to the probabilities of completing the
corresponding pattern at \emph{one} of these trials (in the case of \textsf{%
SSFFS}, only Trials $n$ and $n-4$ are eligible), thus:%
\begin{equation}
p^{3}q^{2}=u_{n}+u_{n-4}\cdot p^{2}q^{2}  \label{Un}
\end{equation}%
which holds for any $n\geq 5$ ($n\geq m$ in general). Multiplying by $s^{n}$
and adding over $n$ from $m$ to infinity (note that $u_{1}$ to $u_{m-1}$
must equal to $0$), yields 
\begin{equation}
\frac{p^{3}q^{2}s^{5}}{1-s}=(1+p^{2}q^{2}s^{4})\left( U(s)\overset{}{-}%
1\right)   \label{U}
\end{equation}%
which can be easily solved for $U(s)-1,$ and then converted to $F(s)$ using (%
\ref{FU}), getting 
\begin{equation*}
F(s)=\frac{1}{1+(1-s)\cdot \dfrac{1+p^{2}q^{2}s^{4}}{p^{3}q^{2}s^{5}}}
\end{equation*}%
Note that $U(1)$ is infinite but $F(1)=1,$ as expected. The mean number of
trials needed to generate the first occurrence of this pattern is obtained
from 
\begin{equation}
\mu \equiv F^{\prime }(s=1)=\frac{1+p^{2}q^{2}}{p^{3}q^{2}}  \label{mu}
\end{equation}

It is difficult to spell out the general form of (\ref{Un}) which would
apply to any pattern, but the idea is (hopefully) quite clear. One should
note that sometimes there is only the $u_{n}$ term on the RHS of the
equation (consider the \textsf{SSFF} pattern), sometimes we have all $m$
terms (e.g. \textsf{SSSS} yields $u_{n}+u_{n-1}\cdot p+u_{n-2}\cdot
p^{2}+u_{n-3}\cdot p^{3}$). The best way to do this is to slide the pattern
past itself to see how many perfect matches one gets (each one of these
contributes exactly one term to the RHS).

\subsection{with a head-start}

To play two or more such patterns against each other (i.e. observing which
of them happens first), we must also find the PGF of the \emph{remaining}
number of trials needed to generate the first occurrence of a pattern, given
that the first few of its symbols are already there. Thus (using the old 
\textsf{SSFFS} example), $F^{\text{SSF}}(s)$ is a PGF of the number of
trials to generate \textsf{SSFFS}, assuming an \textsf{SSF} head start. This
implies that, if we are lucky, we can complete \textsf{SSFFS} in only $2$
more trials (i.e. now $f_{2}^{\text{SSF}}=qp$), and that is pretty much the
only help we can get from \textsf{SSF} in this case. Note that this may get
more complicated in general - we may get another, `shorter' help from the
head-start string (consider \textsf{FSSFF} with the head start of \textsf{%
FSSF} - now we can finish the pattern in a single trial by getting an 
\textsf{F}, but if we get an \textsf{S} instead, we can still complete the
pattern in only three more trials).

Using a similar approach to deriving (\ref{UF}), the corresponding formula
now reads%
\begin{equation}
u_{n}^{\text{SSF}}=\sum_{i=0}^{n}f_{i}^{\text{SSF}}u_{n-i}  \label{xxx}
\end{equation}%
valid for \emph{all} $n\geq 0\,$, since now we take $u_{0}^{\text{SSF}}=0$
(this applies to all $u_{0}$ with a superscript; $u_{0}=1$ remains an
exception). Multiplying (\ref{xxx}) by $s^{n}$ and summing over all $n$
(this time, we \emph{include} $n=0$) yields%
\begin{equation}
F^{\text{SSF}}(s)=\frac{U^{\text{SSF}}(s)}{U(s)}  \label{yyy}
\end{equation}%
where $U(s)$ was defined in the previous section. To get $U^{\text{SSF}}(s),$
we repeat the logic of (\ref{Un}), getting an identical%
\begin{equation}
p^{3}q^{2}=u_{n}^{\text{SSF}}+u_{n-4}^{\text{SSF}}\cdot p^{2}q^{2}
\label{zzz}
\end{equation}%
for $n\geq m$. What changes is that, instead of the old $u_{0}=1,$ $%
u_{1}=u_{2}=u_{3}=u_{4}=0,$ we now have $u_{0}^{\text{SSF}}=u_{1}^{\text{SSF}%
}=u_{3}^{\text{SSF}}=u_{4}^{\text{SSF}}=0$ but $u_{2}^{\text{SSF}}=pq.$
Multiplying (\ref{zzz}) by $s^{n}$ and summing over $n$ from $m$ to infinity
yields%
\begin{equation*}
\frac{p^{3}q^{2}s^{5}}{1-s}=\left( U^{\text{SSF}}(s)\overset{}{-}%
pqs^{2}\right) +p^{2}q^{2}s^{4}\cdot U^{\text{SSF}}(s)
\end{equation*}%
which can be easily solved for $U^{\text{SSF}}(s)$ and consequently
converted to 
\begin{equation*}
F^{\text{SSF}}(s)=\frac{(1-s+p^{2}qs^{3})pqs^{2}}{%
1-s+p^{2}q^{2}s^{4}-p^{2}q^{3}s^{5}}
\end{equation*}%
$\allowbreak $after some simplification. Note that, similarly to $F(s)$, the
new $F^{\text{SSF}}(s)$ also evaluates to $1$ at $s=1.$ This is a universal
property of these PGFs, indicating the any pattern will be generated with
the probability of $1$ sooner or later.

The corresponding mean number of flips is now slightly smaller than (\ref{mu}%
), namely%
\begin{equation*}
\mu ^{\text{SSF}}=\left. \frac{d}{ds}F^{\text{SSF}}(s)\right\vert _{s=1}=%
\frac{1+p^{2}q^{2}-pq}{p^{3}q^{2}}
\end{equation*}

\subsubsection{Another example}

To find $U^{\text{SSS}}(s)$ for the \textsf{SSSS} pattern we start with%
\begin{equation*}
p^{5}=u_{n}^{\text{SSS}}+u_{n-1}^{\text{SSS}}\cdot p+u_{n-2}^{\text{SSS}%
}\cdot p^{2}+u_{n-3}^{\text{SSS}}\cdot p^{3}
\end{equation*}%
valid for $n\geq 4,$ realize that $u_{0}^{\text{SSS}}=u_{2}^{\text{SSS}%
}=u_{3}^{\text{SSS}}=0$, $u_{1}^{\text{SSS}}=p$, and end up with (after
multiplying the previous equation by $s^{n}$ and summing over $n$ from $4$
to infinity):%
\begin{equation*}
\frac{p^{4}s^{4}}{1-s}=\left( U^{\text{SSS}}(s)\overset{}{-}ps\right)
+ps\cdot \left( U^{\text{SSS}}(s)\overset{}{-}sp\right) +p^{2}s^{2}\cdot
\left( U^{\text{SSS}}(s)\overset{}{-}sp\right) +p^{3}s^{3}\cdot U^{\text{SSSS%
}}(s)
\end{equation*}%
which can be easily solved for $U^{\text{SSS}}(s)$ and converted to $F^{%
\text{SSS}}(s).$

\section{Playing 2 patterns against each other}

\subsection{from scratch}

Let us consider two patterns which may not necessarily be of the same
length, but neither of them is allowed to be a substring of the other.

We now modify our notation: let $F_{1}(s)$ and $F_{2}(s)$ will be the PGFs
of the number of trials to generate Pattern $1$ and Pattern $2$
(respectively) for the first time \emph{from scratch}, while $F_{1|2}(s)$
assumes that Pattern $2$ has just been completed and can be used as a head
start to help generate Pattern $1$; similarly, we define $F_{2|1}(s)$. The
corresponding expected values will be denoted $\mu _{1},$ $\mu _{2},$ $\mu
_{1|2}$ and $\mu _{2|1}$ respectively.

Thus, for example, if the first pattern is \textsf{SSFFS} and the second one
is \textsf{FSFSSF, }$F_{1}(s)$ and $F_{1|2}(s)$ are the same as $F(s)\ $and $%
F^{\text{SSF}}(s)$ of the previous section, since the first pattern can use
only the last three symbols of the second pattern (namely \textsf{SSF}) as
its head start. Similarly, the second pattern can use only the last \textsf{%
two} symbols of the first pattern, which then defines $F_{2|1}(s)$; finding
it, together with $F_{2}(s),$ is left as an exercise. Here we quote only the
corresponding two means: 
\begin{eqnarray*}
\mu _{2} &=&\frac{1+p^{3}q^{2}}{p^{3}q^{3}} \\
\mu _{2|1} &=&\frac{1-p^{2}q^{3}}{p^{3}q^{3}}
\end{eqnarray*}

Let now $X_{1\{2\}}(s)$ be the SGF of the probabilities that Pattern $1$
wins over Pattern $2$ at the completion of the $n^{\text{th}}$ trial, and
let $X_{2\{1\}}(s)$ be its vice-versa counterpart. Clearly, $%
x_{1\{2\},0}=x_{2\{1\},0}=0.$

The $f_{1,n}$ probability (of Pattern $1$ completed, for the first time, at
Trial $n$) can be expanded as follows:%
\begin{equation}
f_{1,n}=x_{1\{2\},n}+\sum_{i=0}^{n}x_{2\{1\},i}\cdot f_{1|2,n-i}  \label{key}
\end{equation}%
which is now correct for \emph{all} $n\geq 0.$ The RHS applies the total
probability formula to the sample space of the first $n$ trials,
partitioning it according the trial (denoted $i$) at which Pattern $2$ wins
the game, adding the probability (the first term of the RHS) that Pattern $2$
has not been completed yet, in which case Pattern $1$ has won, at Trial $n$.
Multiplying the previous equation by $s^{n}$ and summing over $n$ (from $0$
to infinity) yields%
\begin{equation*}
F_{1}(s)=X_{1\{2\}}(s)+X_{2\{1\}}(s)\cdot F_{1|2}(s)
\end{equation*}%
since the last term of (\ref{key}) is a \emph{convolution} of the two
sequences, becoming a \emph{product} of the corresponding generating
functions, as explained earlier.

The same must be true with Patterns $1$ and $2$ interchanged, thus:%
\begin{equation*}
F_{2}(s)=X_{2\{1\}}(s)+X_{1\{2\}}(s)\cdot F_{2|1}(s)
\end{equation*}%
The last two equations are easily solved for%
\begin{eqnarray}
X_{1\{2\}}(s) &=&\frac{F_{1}(s)-F_{2}(s)\cdot F_{1|2}(s)}{1-F_{1|2}(s)\cdot
F_{2|1}(s)}  \label{zero} \\
X_{2\{1\}}(s) &=&\frac{F_{2}(s)-F_{1}(s)\cdot F_{2|1}(s)}{1-F_{1|2}(s)\cdot
F_{2|1}(s)}  \notag
\end{eqnarray}

The probability that Pattern $1$ wins (at some trial) is given by $%
X_{1\{2\}}(1),$ or more accurately (since a simple evaluation would lead to
an indefinite answer of $\frac{0}{0}$) by 
\begin{equation*}
\lim_{s\rightarrow 1}X_{1\{2\}}(s)=\frac{\mu _{2}-\mu _{1}+\mu _{1|2}}{\mu
_{1|2}+\mu _{2|1}}
\end{equation*}%
Applied to our example of playing \textsf{SSFFS} against \textsf{FSFSSF}
this yields%
\begin{equation*}
\frac{1-pq^{3}(1+p)}{1+q^{2}+p^{2}q}
\end{equation*}%
which increases from $50\%$ in the $p\rightarrow 0$ limit to $100\%$ when $%
p\rightarrow 1$ (note that Pattern $1$ consists of the same number of 
\textsf{S}s as Pattern $2,$ but has fewer \textsf{F}s).

\subsection{with a head start}

To get ready for playing \emph{three} patterns against each other (the main
topic of this article), we have to extend the previous two formulas by
assuming that the game (of playing Pattern $1$ against Pattern $2$) started
from a completed Pattern $3$ (and allowing either of the two competing
patterns to utilize any of its symbols). By the same reasoning, it is easy
to find that%
\begin{eqnarray}
X_{1\{2\}|3}(s) &=&\frac{F_{1|3}(s)-F_{2|3}(s)\cdot F_{1|2}(s)}{%
1-F_{1|2}(s)\cdot F_{2|1}(s)}  \label{one} \\
X_{2\{1\}|3}(s) &=&\frac{F_{2|3}(s)-F_{1|3}(s)\cdot F_{2|1}(s)}{%
1-F_{1|2}(s)\cdot F_{2|1}(s)}  \notag
\end{eqnarray}%
where $X_{1\{2\}|3}(s)$ is the SGF of the probabilities of Pattern $1$
winning over Pattern $2$ in exactly $n$ trials, \emph{given} that the game
has started from a completed Pattern $3.$

\section{Playing 3 patterns}

We can extend (\ref{key}) by expanding the probability of Pattern $1$
beating Pattern $2$ at a completion of Trial $n$ (the LHS), partitioning the
sample space according the trial (denoted $i$) at which Pattern $3$ has been
completed \emph{for the first time}, including the possibility (the first
term on the RHS), that Pattern $3$ has not been completed yet. In a similar
manner to deriving (\ref{key}), this leads to%
\begin{equation}
X_{1\{2\}}(s)=X_{1\{2,3\}}(s)+X_{3\{2,1\}}(s)\cdot X_{1\{2\}|3}(s)
\label{main}
\end{equation}%
where $X_{1\{2,3\}}(s)$ is the SGF of the probability of Pattern $1$ beating
both Patterns $2$ \emph{and} $3$ at the completion of the $n^{\text{th}}$
trial.

Reversing the r\^{o}le of Patterns $1$ and $3$ we get%
\begin{equation*}
X_{3\{2\}}(s)=X_{3\{2,1\}}(s)+X_{1\{2,3\}}(s)\cdot X_{3\{2\}|1}(s)
\end{equation*}%
Solving the last two equations, one gets 
\begin{eqnarray*}
X_{1\{2,3\}}(s) &=&\frac{X_{1\{2\}}(s)-X_{3\{2\}}(s)\cdot X_{1\{2\}|3}(s)}{%
1-X_{3\{2\}|1}(s)\cdot X_{1\{2\}|3}(s)} \\
X_{3\{2,1\}}(s) &=&\frac{X_{3\{2\}}(s)-X_{1\{2\}}(s)\cdot X_{3\{2\}|1}(s)}{%
1-X_{1\{2\}|3}(s)\cdot X_{3\{2\}|1}(s)}
\end{eqnarray*}%
and, by permuting the indices, an analogous solution for $X_{1\{3,2\}},$ $%
X_{3\{1,2\}},$ $X_{2\{1,3\}}$ and $X_{2\{3,1\}}$ (at this point it would
appear that $X_{1\{2,3\}}(s)$ is different from $X_{1\{3,2\}}(s),$ but keep
on reading).

Utilizing (\ref{zero}) and (\ref{one}), we can then express any of these
solutions directly in terms of the $F_{i}\left( s\right) $ and $F_{i|j}(s)$
functions, getting (for expedience, we quote each $F$ without its $(s)$
argument):%
\begin{eqnarray}
&&X_{1\{2,3\}}(s)\overset{}{=}  \label{sol} \\
&&\frac{F_{1}\left( 1-F_{2|3}F_{3|2}\right) -F_{2}\left(
F_{1|2}-F_{1|3}F_{3|2}\right) -F_{3}\left( F_{1|3}-F_{1|2}F_{2|3}\right) }{%
1-F_{1|2}F_{2|1}-F_{1|3}F_{3|1}-F_{2|3}F_{3|2}+F_{1|2}F_{2|3}F_{3|1}+F_{1|3}F_{3|2}F_{2|1}%
}  \notag
\end{eqnarray}%
and its index-permuted equivalents. Now, it becomes explicitly obvious that $%
X_{1\{2,3\}}(s)=X_{1\{3,2\}}(s)$ as expected (winning over Patterns $2$ and $%
3$ is the same as winning over Patterns $3$ and $2$).

To get the probability of Pattern $1$ winning over the other two patterns at 
\emph{any} trial, one has to evaluate (\ref{sol}) at $s=1.$ This yields
(with the help of L'Hospital rule, having to differentiate each numerator
and denominator \emph{twice}) 
\begin{eqnarray}
&&\lim_{s\rightarrow 1}X_{1\{2,3\}}(s)\overset{}{=}  \label{ED} \\
&&\frac{%
\begin{array}{c}
\mu _{1}(\mu _{2|3}+\mu _{3|2})+\mu _{2}(\mu _{1|2}-\mu _{1|3}-\mu
_{3|2})+\mu _{3}(\mu _{1|3}-\mu _{1|2}-\mu _{2|3}) \\ 
+\mu _{2|3}\mu _{3|2}-\mu _{1|3}\mu _{3|2}-\mu _{1|2}\mu _{2|3}%
\end{array}%
}{%
\begin{array}{c}
\mu _{1|2}\mu _{2|1}+\mu _{1|3}\mu _{3|1}+\mu _{2|3}\mu _{3|2}-\mu _{1|2}\mu
_{2|3}-\mu _{1|3}\mu _{3|2} \\ 
-\mu _{2|1}\mu _{1|3}-\mu _{2|3}\mu _{3|1}-\mu _{3|1}\mu _{1|2}-\mu
_{3|2}\mu _{2|1}%
\end{array}%
}  \notag
\end{eqnarray}%
One can get the other two answers by permuting indices.

Finally, the PGF of the game's duration is clearly given by the following sum%
\begin{equation*}
H(s)\equiv X_{(1|2,3)}(s)+X_{(2|1,3)}(s)+X_{(3|1,2)}(s)
\end{equation*}

The expected duration of the game is thus equal to the $s$ derivative of
this expression, evaluated (with the help of L'Hospital rule and the \emph{%
fourth} derivative of the numerator and denominator) at $s=1.$ This yields 
\begin{eqnarray}
&&\lim_{s\rightarrow 1}H^{\prime }(s)\overset{}{=}  \label{dur} \\
&&\frac{%
\begin{array}{c}
\mu _{1}(\mu _{2|3}\mu _{3|2}-\mu _{2|3}\mu _{3|1}-\mu _{3|2}\mu _{2|1})+\mu
_{2}(\mu _{1|3}\mu _{3|1}-\mu _{1|3}\mu _{3|2}-\mu _{3|1}\mu _{1|2}) \\ 
+\mu _{3}(\mu _{1|2}\mu _{2|1}-\mu _{1|2}\mu _{2|3}-\mu _{2|1}\mu
_{1|3})+\mu _{1|2}\mu _{2|3}\mu _{3|1}+\mu _{1|3}\mu _{3|2}\mu _{2|1}%
\end{array}%
}{%
\begin{array}{c}
\mu _{1|2}\mu _{2|1}+\mu _{1|3}\mu _{3|1}+\mu _{2|3}\mu _{3|2}-\mu _{1|2}\mu
_{2|3}-\mu _{1|3}\mu _{3|2} \\ 
-\mu _{2|1}\mu _{1|3}-\mu _{2|3}\mu _{3|1}-\mu _{3|1}\mu _{1|2}-\mu
_{3|2}\mu _{2|1}%
\end{array}%
}  \notag
\end{eqnarray}

\subsection{Example}

Adding \textsf{FSSSF} to the two patterns of the previous section, we compute%
\begin{eqnarray*}
\mu _{3} &=&\frac{1+p^{3}q}{p^{3}q^{2}} \\
\mu _{3|1} &=&\mu _{3}^{\text{FS}}=\frac{1-p^{2}q^{2}}{p^{3}q^{2}} \\
\mu _{3|2} &=&\mu _{3}^{\text{F}}=\frac{1}{p^{3}q^{2}} \\
\mu _{1|3} &=&\mu _{1}^{\text{SSF}}=\mu _{1|2} \\
\mu _{2|3} &=&\mu _{2}^{\text{F}}=\frac{1}{p^{3}q^{3}}
\end{eqnarray*}%
Based on (\ref{ED}), the probability of \textsf{SSFFS} winning the game is%
\begin{equation*}
\frac{1-pq^{2}(1+p)(1+q)}{3q+p^{2}(2+q)}
\end{equation*}%
which varies from $\frac{1}{3}$ at $p\rightarrow 0$ (all three patterns have
the same chance of winning, since they contain the same number of \textsf{S}%
s, and that is all what counts in this limit) to its smallest value of $%
28.59\%$ at $p=0.2495,$ to its largest value of $50\%$ at $p\rightarrow 1$
(Patterns $1$ and $3$ have the same chance of winning, each containing two 
\textsf{F}s; Pattern $2$ with three \textsf{F}s is now out of the contest).

The expected duration of the game is%
\begin{equation*}
\frac{1+p^{2}q\left( 1\overset{}{-}pq^{3}(1+pq)\right) }{p^{3}q^{2}\left( 3q%
\overset{}{+}p^{2}(2+q)\right) }
\end{equation*}%
which reaches its smallest value of $15.88$ trials at $p=0.5796$; in each of
the $p\rightarrow 0$ and $p\rightarrow 1$ limits, this mean (or average)
number of trials to complete this game becomes infinite.

\subsection{Final challenge}

It would be an interesting (but certainly non-trivial) exercise to extend
the theory to four or more competing patterns.

\end{document}